\documentclass[11pt,a4paper,twoside]{article}
\usepackage{latexsym,amsmath,amsthm,amssymb,amsfonts}
\usepackage{mathrsfs,multicol}
\usepackage{amscd}
\usepackage{cite}
\usepackage[ansinew]{inputenc}

\def\proof{{\bf Proof.}\quad}

\numberwithin{equation}{section}

\catcode`@=11

\newskip\plaincentering \plaincentering=0pt plus 1000pt minus 1000pt
\def\@plainlign{\tabskip=0pt\everycr={}}
\def\eqalignno#1{\displ@y \tabskip\plaincentering
  \halign to\displaywidth{\hfil$\@lign\displaystyle{##}$\tabskip\z@skip
    &$\@lign\displaystyle{{}##}$\hfil\tabskip\plaincentering
    &\llap{$\@lign##$}\tabskip\z@skip\crcr
    #1\crcr}}
\def\leqalignno#1{\displ@y \tabskip\plaincentering
  \halign to\displaywidth{\hfil$\@lign\displaystyle{##}$\tabskip\z@skip
    &$\@lign\displaystyle{{}##}$\hfil\tabskip\plaincentering
    &\kern-\displaywidth\rlap{$\@lign##$}\tabskip\displaywidth\crcr
    #1\crcr}}
\def\plainLet@{\relax\iffalse{\fi\let\\=\cr\iffalse}\fi}
\def\plainvspace@{\def\vspace##1{\noalign{\vskip##1}}}

\def\intic@{\mathchoice{\hskip5\p@}{\hskip4\p@}{\hskip4\p@}{\hskip4\p@}}
\def\negintic@
 {\mathchoice{\hskip-5\p@}{\hskip-4\p@}{\hskip-4\p@}{\hskip-4\p@}}
\def\intkern@{\mathchoice{\!\!\!}{\!\!}{\!\!}{\!\!}}
\def\intdots@{\mathchoice{\cdots}{{\cdotp}\mkern1.5mu
    {\cdotp}\mkern1.5mu{\cdotp}}{{\cdotp}\mkern1mu{\cdotp}\mkern1mu
      {\cdotp}}{{\cdotp}\mkern1mu{\cdotp}\mkern1mu{\cdotp}}}
\newcount\intno@
\def\iint{\intno@=\tw@\futurelet\next\ints@}
\def\iiint{\intno@=\thr@@\futurelet\next\ints@}
\def\iiiint{\intno@=4 \futurelet\next\ints@}
\def\idotsint{\intno@=\z@\futurelet\next\ints@}
\def\ints@{\findlimits@\ints@@}
\newif\iflimtoken@
\newif\iflimits@
\def\findlimits@{\limtoken@false\limits@false\ifx\next\limits
 \limtoken@true\limits@true\else\ifx\next\nolimits\limtoken@true\limits@false
    \fi\fi}
\def\multintlimits@{\intop\ifnum\intno@=\z@\intdots@
  \else\intkern@\fi
    \ifnum\intno@>\tw@\intop\intkern@\fi
     \ifnum\intno@>\thr@@\intop\intkern@\fi\intop}
\def\multint@{\int\ifnum\intno@=\z@\intdots@\else\intkern@\fi
   \ifnum\intno@>\tw@\int\intkern@\fi
    \ifnum\intno@>\thr@@\int\intkern@\fi\int}
\def\ints@@{\iflimtoken@\def\ints@@@{\iflimits@
   \negintic@\mathop{\intic@\multintlimits@}\limits\else
    \multint@\nolimits\fi\eat@}\else
     \def\ints@@@{\multint@\nolimits}\fi\ints@@@}
\def\Sb{_\bgroup\vspace@
        \baselineskip=\fontdimen10 \scriptfont\tw@
        \advance\baselineskip by \fontdimen12 \scriptfont\tw@
        \lineskip=\thr@@\fontdimen8 \scriptfont\thr@@
        \lineskiplimit=\thr@@\fontdimen8 \scriptfont\thr@@
        \Let@\vbox\bgroup\halign\bgroup \hfil$\scriptstyle
            {##}$\hfil\cr}
\def\endSb{\crcr\egroup\egroup\egroup}
\def\Sp{^\bgroup\vspace@
        \baselineskip=\fontdimen10 \scriptfont\tw@
        \advance\baselineskip by \fontdimen12 \scriptfont\tw@
        \lineskip=\thr@@\fontdimen8 \scriptfont\thr@@
        \lineskiplimit=\thr@@\fontdimen8 \scriptfont\thr@@
        \Let@\vbox\bgroup\halign\bgroup \hfil$\scriptstyle
            {##}$\hfil\cr}
\def\endSp{\crcr\egroup\egroup\egroup}
\def\Let@{\relax\iffalse{\fi\let\\=\cr\iffalse}\fi}
\def\vspace@{\def\vspace##1{\noalign{\vskip##1 }}}
\def\aligned{\,\vcenter\bgroup\plainvspace@\plainLet@\openup\jot\m@th\ialign
  \bgroup \strut\hfil$\displaystyle{##}$&$\displaystyle{{}##}$\hfil\crcr}
\def\endaligned{\crcr\egroup\egroup}
\def\matrix{\,\vcenter\bgroup\plainLet@\plainvspace@
    \normalbaselines
  \m@th\ialign\bgroup\hfil$##$\hfil&&\quad\hfil$##$\hfil\crcr
    \mathstrut\crcr\noalign{\kern-\baselineskip}}
\def\endmatrix{\crcr\mathstrut\crcr\noalign{\kern-\baselineskip}\egroup
                \egroup\,}
\newtoks\hashtoks@
\hashtoks@={#}
\def\format{\crcr\egroup\iffalse{\fi\ifnum`}=0 \fi\format@}
\def\format@#1\\{\def\preamble@{#1}%
  \def\c{\hfil$\the\hashtoks@$\hfil}%
  \def\r{\hfil$\the\hashtoks@$}%
  \def\l{$\the\hashtoks@$\hfil}%
  \setbox\z@=\hbox{\xdef\Preamble@{\preamble@}}\ifnum`{=0 \fi\iffalse}\fi
   \ialign\bgroup\span\Preamble@\crcr}

\def\cases{\left\{\,\vcenter\bgroup\plainvspace@
     \normalbaselines\openup\jot\m@th
      \plainLet@\ialign\bgroup$\displaystyle{##}$\hfil&
      \quad$\displaystyle{{}##}$\hfil\crcr
      \mathstrut\crcr\noalign{\kern-\baselineskip}}

\newif\iftagsleft@
\tagsleft@true
\def\TagsOnRight{\global\tagsleft@false}
\def\tag#1$${\iftagsleft@\leqno\else\eqno\fi
 \hbox{\def\pagebreak{\global\postdisplaypenalty-\@M}%
 \def\nopagebreak{\global\postdisplaypenalty\@M}\rm(#1\unskip)}%
  $$\postdisplaypenalty\z@\ignorespaces}
\interdisplaylinepenalty=\@M
\def\allowdisplaybreak{\noalign{\allowbreak}}
\def\plainallowdisplaybreak@{\def\allowdisplaybreak{\noalign{\allowbreak}}}
\def\plaindisplaybreak@{\def\displaybreak{\noalign{\break}}}
\def\align#1\endalign{\def\tag{&}\plainvspace@\plainallowdisplaybreak@
\plaindisplaybreak@
  \iftagsleft@\plainlalign@#1\endalign\else
   \plainralign@#1\endalign\fi}
\def\plainralign@#1\endalign{\displ@y\plainLet@\tabskip\plaincentering
\halign to\displaywidth
     {\hfil$\displaystyle{##}$\tabskip=\z@&$\displaystyle{{}##}$\hfil
       \tabskip=\plaincentering&\llap{\hbox{\rm(##\unskip)}}\tabskip\z@\crcr
             #1\crcr}}
\def\plainlalign@
 #1\endalign{\displ@y\plainLet@\tabskip\plaincentering\halign to \displaywidth
   {\hfil$\displaystyle{##}$\tabskip=\z@&$\displaystyle{{}##}$\hfil
   \tabskip=\plaincentering&\kern-\displaywidth
        \rlap{\hbox{\rm(##\unskip)}}\tabskip=\displaywidth\crcr
               #1\crcr}}

\def\re@#1{\par\hangindent\parindent\indent\llap{#1\enspace}\ignorespaces}
\def\qfootnote#1{\edef\@sf{\spacefactor\the\spacefactor}{}#1\@sf
      \insert\footins{\let\egroup=}\footnotesize 
      \interlinepenalty100 \let\par=\endgraf
        \leftskip=0pt \rightskip=0pt
        \splittopskip=10pt plus 1pt minus 1pt \floatingpenalty=20000
   \smallskip\re@{#1}\bgroup\strut\aftergroup{\strut\egroup}\let\next}
\catcode`\@=\active
\TagsOnRight
\begin{document}
\title{\bf On a classification of the quasi Yamabe gradient solitons
\footnote{The research of the first author is supported by NSFC
grant(No. 11001076) and NSF of Henan Provincial Education
department(No. 2010A110008). The research of the second author is
supported by NSFC grant(No. 10971110).}}
\author{Guangyue Huang and Haizhong Li
\footnote{The corresponding author's
email:~\textsf{hli$@$math.tsinghua.edu.cn}\,(H. Li)}}
\date{}
\maketitle
\begin{quotation}
\noindent{\bf Abstract.}~ In this paper, we introduce the concept of quasi Yamabe gradient solitons,
which generalizes the concept of Yamabe gradient solitons. By using some ideas in \cite{CaoChen1,Caoetal1},
we prove that $n$-dimensional $(n\geq3)$
complete quasi Yamabe gradient solitons
with vanishing Weyl curvature tensor and positive sectional curvature must be rotationally symmetric. We also
prove that any compact quasi Yamabe gradient solitons
are of constant scalar curvature.
\end{quotation}

\section{Introduction}

Let $(M^n,g)$ be an $n$-dimensional Riemannian manifold with
$n\geq3$. If there exists a smooth function $f$ on $M^n$ and a
constant $\rho$ such that
\begin{equation}\label{Int1}
(R-\rho)g_{ij}=f_{ij},
\end{equation} then we call $(M^n,g,f)$ a {\it Yamabe gradient soliton}.
Here $R$ denotes the scalar curvature of the metric $g$. If
$\rho=0$, $\rho>0$ or $\rho<0$, then $(M^n,g,f)$ is called a Yamabe
steady, Yamabe shrinker or Yamabe expander respectively. Yamabe
solitons are special solutions to the Yamabe flow
\begin{equation}\label{Int2}
\frac{\partial}{\partial t} g_{ij}=-Rg_{ij}.
\end{equation}
For the study of the Yamabe flow in the compact case, see
\cite{Chow92,Brendle2007,Brendle2005,Ye1994,delPino01,Schwetlick03}
and the references therein. It is very important for understanding
the singularity formation in the complete Yamabe flow to study the
classification of Yamabe solitons.

In \cite{Sesum11}, Daskalopoulos and Sesum studied the
classification of locally conformally flat Yamabe gradient solitons.
They proved

\noindent{\bf Theorem A.(\cite{Sesum11})} {\it If $(M,g,f)$ is a
complete locally conformally flat Yamabe gradient soliton satisfying
\eqref{Int1} with positive sectional curvature, then $(M,g,f)$ is
rotationally symmetric.

}

\noindent{\bf Theorem B.(\cite{Sesum11})} {\it If $(M,g,f)$ is a
compact Yamabe gradient soliton, then $g$ is the metric of constant
scalar curvature.

}

In this paper, we consider a generalized Yamabe gradient soliton which
we call the {\it quasi Yamabe gradient soliton}.

\noindent{\bf Definition 1.1.} {\it If there exists a smooth function
$f$ on $M^n$ and two constants $m,\rho$ (where $m$ is not zero) such
that
\begin{equation}\label{Int3}
(R-\rho)g_{ij}=f_{ij}-\frac{1}{m}f_if_j,
\end{equation}
then we call $(M^n,g,f)$ a quasi Yamabe gradient soliton.

}

We remark that if $m\rightarrow\infty$,  \eqref{Int3} reduces to \eqref{Int1}, so quasi Yamabe gradient solitons can be considered as generalized Yamabe gradient solitons in this sense.

We study classifications for complete quasi Yamabe gradient solitons
satisfying \eqref{Int3}. As in \cite{CaoChen1,Caoetal1}, the key idea in proving our results is to
link the Weyl curvature tensor with the covariant 3-tensor $D_{ijk}$, introduced by Cao-Chen \cite{CaoChen1,CaoChen2}, with
$D_{ijk}=W_{ijkl}f_l$, see Proposition 2.2 in Section 2, where
$D_{ijk}$ is defined by
\begin{equation}\label{Int4}\aligned
D_{ijk}=&\frac{1}{n-2}(R_{kj}f_i-R_{ki}f_j)+\frac{1}{(n-1)(n-2)}
(R_{il}g_{jk}f^l-R_{jl}g_{ik}f^l)\\
&-\frac{R}{(n-1)(n-2)}(g_{kj}f_i-g_{ki} f_j).
\endaligned\end{equation}
Our main results are as follows:

\noindent{\bf Theorem 1.1.} {\it Let $(M,g,f)$ be a complete quasi
Yamabe gradient soliton satisfying \eqref{Int3} with positive
sectional curvature.

(1) If $n=3$, then $(M,g,f)$ is rotationally symmetric;

(2) If $n=4$ and $D_{ijk}=0$, then $(M,g,f)$ is rotationally
symmetric;

(3) If $n\geq5$ and $W_{ijkl}=0$, then $(M,g,f)$ is rotationally
symmetric.

}

\noindent{\bf Theorem 1.2.} {\it If $(M,g,f)$ is a compact quasi
Yamabe gradient soliton satisfying \eqref{Int3}, then $g$ is the
metric of constant scalar curvature.

}

Letting $m\rightarrow \infty$, we have the following result from Theorem 1.1 immediately:

\noindent{\bf Corollary 1.1.} {\it Let $(M,g,f)$ be a complete
Yamabe gradient soliton satisfying \eqref{Int1} with positive
sectional curvature.

(1) If $n=3$, then $(M,g,f)$ is rotationally symmetric;

(2) If $n=4$ and $D_{ijk}=0$, then $(M,g,f)$ is rotationally
symmetric;

(3) If $n\geq5$ and $W_{ijkl}=0$, then $(M,g,f)$ is rotationally
symmetric.

}

\noindent{\it Remark 1.1.} Daskalopoulos and Sesum in \cite{Sesum11}
proved Theorem A (see Theorem 1.3 in \cite{Sesum11}) under the
assumption that the metric $g$ is locally conformally flat. Theorem A follows from Corollary 1.1.
Theorem B follows from Theorem 1.2 when $m\rightarrow \infty$.

\noindent{\it Remark 1.2.} Since the 3-tensor $D_{ijk}$ is related
to the Weyl curvature tensor by $D_{ijk}=W_{ijkl}f^l$ (see
Proposition 2.2 in Section 2), we have from Proposition 2.5 in
Section 2 that $D_{ijk}=0$ is equivalent to $W_{ijkl}f^l=0$ when
$n=4$. However, for $n\geq5$, we can not conclude that $D_{ijk}=0$
implies $W_{ijkl}f^l=0$.

\noindent{\it Remark 1.3.} Some related results for the gradient Ricci solitons can be found in \cite{ Petersen10,Ni08,Nilei08,Chenwang11,Caohuai11,Caohuai10} and the references therein.

\section{Proof of Theorem 1.1}

Throughout this paper, we will agree on the following index
convention:
$$1\leq i,j,k,\cdots\leq n;\ \ \ \ \ \ 2\leq a,b,c,\cdots\leq n.$$
For convenience, we define
\begin{equation}\label{Proof1}
R_\rho=R-\rho.
\end{equation} Then \eqref{Int3} becomes
\begin{equation}\label{Proof2}
R_\rho\,g_{ij}=f_{ij}-\frac{1}{m}f_if_j.
\end{equation}
For quasi Yamabe gradient solitons, we have the following lemma
which will be used later:

\noindent{\bf Lemma 2.1.} {\it Let $(M^n,g,f)$ be a quasi Yamabe
gradient soliton satisfying \eqref{Int3}. Then we have
\begin{equation}\label{Proof3}
nR_\rho=\Delta f-\frac{1}{m}|\nabla f|^2,
\end{equation}
\begin{equation}\label{Proof4}
(|\nabla f|^2)_{i}=2R_\rho f_i+\frac{2}{m}|\nabla f|^2f_i,
\end{equation}
\begin{equation}\label{Proof5}
(R_\rho)_{i}=\frac{1}{m}R_\rho f_i-\frac{1}{n-1}R_{ij}f^j,
\end{equation} where $R_{ij}$ denotes the Ricci curvature of the
metric $g$ and $f^j=g^{jk}f_k$. }

\proof The relationship \eqref{Proof3} can be obtained directly by
contracting the equation \eqref{Proof2}. On the other hand, by choosing the local orthogonal frame $\{e_1,\cdots,e_n\}$, we have by use of
\eqref{Proof2},
$$(|\nabla
f|^2)_{i}=2f_{ij}f_j=2(R_\rho\,g_{ij}+\frac{1}{m}f_if_j)f_j=2R_\rho
f_i+\frac{2}{m}|\nabla f|^2f_i.
$$ Hence, we obtain \eqref{Proof4}.

Using the equation \eqref{Proof2} again, we obtain
\begin{equation}\label{Proof6}
(R_\rho)_{i}=f_{ijj}-\frac{1}{m}f_{ij}f_j-\frac{1}{m}f_if_{jj}.\end{equation}
With the help of the Ricci identity, \eqref{Proof3} and
\eqref{Proof4}, we deduce from \eqref{Proof6}
$$\aligned
(R_\rho)_{i}=&f_{ijj}-\frac{1}{m}f_{ij}f_j-\frac{1}{m}f_if_{jj}\\
=&(\Delta
f)_i+R_{ij}f_j-\frac{1}{2m}(|\nabla f|^2)_{i}
-\frac{1}{m}f_i(\Delta f)\\
=&(nR_\rho+\frac{1}{m}|\nabla f|^2)_i+R_{ij}f_j-\frac{1}{2m}(|\nabla
f|^2)_{i}
-\frac{1}{m}f_i(nR_\rho+\frac{1}{m}|\nabla f|^2)\\
=&n(R_\rho)_{i}+R_{ij}f_j+\frac{1}{2m}(|\nabla f|^2)_{i}
-\frac{n}{m}R_\rho f_i-\frac{1}{m^2}|\nabla f|^2f_i\\
=&n(R_\rho)_{i}+R_{ij}f_j+\frac{1}{2m}(2R_\rho
f_i+\frac{2}{m}|\nabla f|^2f_i)
-\frac{n}{m}R_\rho f_i-\frac{1}{m^2}|\nabla f|^2f_i\\
=&n(R_\rho)_{i}+R_{ij}f_j-\frac{n-1}{m}R_\rho f_i,
\endaligned$$
which concludes the proof of \eqref{Proof5}. It completes the proof
of Lemma 2.1.

For $n\geq3$, the Weyl curvature tensor is defined by
\begin{equation}\label{Proof7}\aligned
W_{ijkl}=&R_{ijkl}-\frac{1}{n-2}(R_{ik}g_{jl}-R_{il}g_{jk}
+R_{jl}g_{ik}-R_{jk}g_{il})\\
+&\frac{R}{(n-1)(n-2)}(g_{ik}g_{jl}-g_{il}g_{jk}).
\endaligned
\end{equation} From the definition of the Weyl curvature tensor
above, it is easy to see that the Weyl curvature tensor satisfies
all the symmetries of the curvature tensor and its traces with the
metric are zero. It is well known that $W_{ijkl}=0$ for $n=3$. For
$n\geq4$, $(M^n,g)$ is locally conformally flat if and only if
$W_{ijkl}=0$.

As in \cite{CaoChen1,CaoChen2}, see also
\cite{CaoChen2,ChenCheng1,Caoetal1,Brendle11}, we define the
following 3-tensor $D$ by
\begin{equation}\label{Proof8}\aligned
D_{ijk}=&\frac{1}{n-2}(R_{kj}f_i-R_{ki}f_j)+\frac{1}{(n-1)(n-2)}
(R_{il}g_{jk}f^l-R_{jl}g_{ik}f^l)\\
&-\frac{R}{(n-1)(n-2)}(g_{kj}f_i-g_{ki} f_j).
\endaligned\end{equation}
Then for quasi Yamabe gradient solitons, we have the following
consequence:

\noindent{\bf Proposition 2.2.} {\it Let $(M^n,g,f)$ be a quasi
Yamabe gradient soliton satisfying \eqref{Int3}. Then the 3-tensor
$D$ is related to the Weyl curvature tensor by
\begin{equation}\label{Proof9}
D_{ijk}=W_{ijkl}f^l.
\end{equation}

}

\proof By use of \eqref{Proof2} and \eqref{Proof5}, we have
\begin{equation}\label{Proof10}\aligned
f_{kji}-f_{kij}=&(R_\rho\,g_{kj}+\frac{1}{m}f_kf_j)_i
-(R_\rho\,g_{ki}+\frac{1}{m}f_kf_i)_j\\
=&(R_\rho)_{i}g_{kj}-(R_\rho)_{j}g_{ki}+\frac{1}{m}(f_{ki}f_j-f_{kj}f_i)\\
=&(R_\rho)_{i}g_{kj}-(R_\rho)_{j}g_{ki}-\frac{R_\rho}{m}(g_{kj}f_i-g_{ki}
f_j)\\
=& -\frac{1}{n-1}(R_{il}g_{jk}f^l-R_{jl}g_{ik}f^l),
\endaligned\end{equation} where the last equality used
\eqref{Proof5}. Using the Ricci identity and \eqref{Proof7}, we have
\begin{equation}\label{Proof11}\aligned
f_{kji}-f_{kij}=&f^lR_{lkji}\\
=&W_{ijkl}f^l-\frac{R}{(n-1)(n-2)}(g_{ik}g_{jl}-g_{il}g_{jk})f^l\\
&+\frac{1}{n-2}(R_{ik}g_{jl}-R_{il}g_{jk}
+R_{jl}g_{ik}-R_{jk}g_{il})f^l\\
=&W_{ijkl}f^l+\frac{R}{(n-1)(n-2)}(g_{kj}f_i-g_{ki} f_j)\\
&-\frac{1}{n-2}(R_{kj}f_i-R_{ki}f_j)-\frac{1}{n-2}
(R_{il}g_{jk}f^l-R_{jl}g_{ik}f^l).
\endaligned\end{equation}
Combining \eqref{Proof10} and \eqref{Proof11}, we obtain by definition \eqref{Proof8} $$\aligned
W_{ijkl}f^l=&\frac{1}{n-2}(R_{kj}f_i-R_{ki}f_j)+\frac{1}{(n-1)(n-2)}
(R_{il}g_{jk}f^l-R_{jl}g_{ik}f^l)\\
&-\frac{R}{(n-1)(n-2)}(g_{kj}f_i-g_{ki} f_j)\\
=&D_{ijk}.
\endaligned$$
Therefore, we complete the proof of Proposition 2.2.

In particular, by properties of the Weyl curvature tensor and
\eqref{Proof8}, we get that $D_{ijk}$ is skew-symmetric in their
first two indices and trace-free in any two indices:
\begin{equation}\label{Proof12}
D_{ijk}=-D_{jik},\ \ \ \ g^{ij}D_{ijk}=g^{ik}D_{ijk}=0.
\end{equation}

Next, we give a proposition which links the norm of $D_{ijk}$ to the
geometry of the level surfaces of the potential function $f$.

\noindent{\bf Proposition 2.3.} {\it Let $(M^n,g,f)$ be a quasi
Yamabe gradient soliton satisfying \eqref{Int3}, and let
$\Sigma_c=\{x|f(x)=c\}$ be the level surface with respect to regular
value $c$ of $f$. Then for any local orthonormal frame $\{e_1, e_2,
\cdots, e_n\}$ with $e_1=\nabla f/|\nabla f|$ and $\{e_2, \cdots,
e_n\}$ tangent to $\Sigma_c$, we have
\begin{equation}\label{Proof13}
|D_{ijk}|^2=\frac{2|\nabla f|^2}{(n-1)(n-2)^2}\Big\{(n-2)R_{1a}^2
+(n-1)\Big|R_{ab}-\frac{R-R_{11}}{n-1}g_{ab}\Big|^2\Big\},
\end{equation}
where $g_{ab}$ is the induced metric on $\Sigma_c$.

}

\proof  From \eqref{Proof8}, we have
\begin{equation}\label{Proof14}\aligned
|D_{ijk}|^2=&\frac{1}{(n-2)^2}|R_{kj}f_i-R_{ki}f_j|^2+\frac{1}{(n-1)^2(n-2)^2}
|R_{il}g_{jk}f_l-R_{jl}g_{ik}f_l|^2\\
&+\frac{R^2}{(n-1)^2(n-2)^2}|g_{kj}f_i-g_{ki}
f_j|^2\\
&+\frac{2}{(n-1)(n-2)^2}(R_{kj}f_i-R_{ki}f_j)(R_{il}g_{jk}f_l-R_{jl}g_{ik}f_l)\\
&-\frac{2R}{(n-1)(n-2)^2}(R_{kj}f_i-R_{ki}f_j)(g_{kj}f_i-g_{ki} f_j)\\
&-\frac{2R}{(n-1)^2(n-2)^2}(R_{il}g_{jk}f_l-R_{jl}g_{ik}f_l)(g_{kj}f_i-g_{ki}
f_j).
\endaligned\end{equation}
Let $\{e_1, e_2, \cdots, e_n\}$ be any local orthonormal frame with
$e_1=\nabla f/|\nabla f|$ and $\{e_2, \cdots, e_n\}$ tangent to
$\Sigma_c$. That is, under this orthonormal frame, we have
$f_1=|\nabla f|$ and $f_2=f_3=\cdots=f_n=0$. Thus,
\begin{equation}\label{Proof15}
|R_{kj}f_i-R_{ki}f_j|^2=2|\nabla f|^2(|Ric|^2-R_{1i}^2),
\end{equation}
\begin{equation}\label{Proof16}
|R_{il}g_{jk}f_l-R_{jl}g_{ik}f_l|^2=2(n-1)|\nabla f|^2R_{1i}^2,
\end{equation}
\begin{equation}\label{Proof17}
|g_{kj}f_i-g_{ki} f_j|^2=2(n-1)|\nabla f|^2,
\end{equation}
\begin{equation}\label{Proof18}
(R_{kj}f_i-R_{ki}f_j)(R_{il}g_{jk}f_l-R_{jl}g_{ik}f_l)=2|\nabla
f|^2(R \,R_{11}-R_{1i}^2),
\end{equation}
\begin{equation}\label{Proof19}
(R_{kj}f_i-R_{ki}f_j)(g_{kj}f_i-g_{ki} f_j)=2|\nabla f|^2(R
-R_{11}),
\end{equation}
\begin{equation}\label{Proof20}
(R_{il}g_{jk}f_l-R_{jl}g_{ik}f_l)(g_{kj}f_i-g_{ki}
f_j)=2(n-1)|\nabla f|^2R_{11}.
\end{equation}
Inserting the relationships \eqref{Proof15}-\eqref{Proof20} into
\eqref{Proof14} yields
\begin{equation}\label{Proof21}\aligned
|D_{ijk}|^2=&\frac{2}{(n-2)^2}|\nabla f|^2(|Ric|^2-R_{1i}^2)
+\frac{2}{(n-1)(n-2)^2}|\nabla f|^2R_{1i}^2\\
&+\frac{2R^2}{(n-1)(n-2)^2}|\nabla f|^2+\frac{4}{(n-1)(n-2)^2}
|\nabla f|^2(R \,R_{11}-R_{1i}^2)\\
&-\frac{4R}{(n-1)(n-2)^2}|\nabla f|^2(R
-R_{11})-\frac{4R}{(n-1)(n-2)^2}|\nabla f|^2R_{11}\\
=&\frac{2|\nabla
f|^2}{(n-1)(n-2)^2}\Big\{(n-1)(|Ric|^2-R_{1i}^2)+R_{1i}^2+R^2+2(R
\,R_{11}-R_{1i}^2)\\
&-2R(R -R_{11})-2R\,R_{11}\Big\}\\
=&\frac{2|\nabla f|^2}{(n-1)(n-2)^2}\Big\{(n-1)|Ric|^2
-nR_{1i}^2+2R\,R_{11}-R^2\Big\}\\
=&\frac{2|\nabla
f|^2}{(n-1)(n-2)^2}\Big\{(n-1)(R_{11}^2+2R_{1a}^2+R_{ab}^2)
-n(R_{11}^2+R_{1a}^2)\\
& +2(R_{11}+R_{aa})R_{11}
-(R_{11}^2+2R_{11}R_{aa}+R_{aa}R_{bb})\Big\}\\
=&\frac{2|\nabla f|^2}{(n-1)(n-2)^2}\Big\{(n-2)R_{1a}^2
+(n-1)\Big|R_{ab}-\frac{R_{cc}}{n-1}g_{ab}\Big|^2\Big\}\\
=&\frac{2|\nabla f|^2}{(n-1)(n-2)^2}\Big\{(n-2)R_{1a}^2
+(n-1)\Big|R_{ab}-\frac{R-R_{11}}{n-1}g_{ab}\Big|^2\Big\}.
\endaligned\end{equation} It completes the proof of Proposition 2.3.

With the help of Proposition 2.3, we can obtain the following
results:

\noindent{\bf Proposition 2.4.} {\it Let $(M^n,g,f)$ be a quasi
Yamabe gradient soliton satisfying \eqref{Int3} with $D_{ijk}=0$,
and let $\Sigma_c=\{x|f(x)=c\}$ be the level surface with respect to
regular value $c$ of $f$. Then for any local orthonormal frame
$\{e_1, e_2, \cdots, e_n\}$ with $e_1=\nabla f/|\nabla f|$ and
$\{e_2, \cdots, e_n\}$ tangent to $\Sigma_c$, we have

(1) $|\nabla f|^2$ and the scalar curvature $R$ of $(M^n, g_{ij},
f)$ are constant on $\Sigma_c$;

(2) $R_{1a}=0$ and $e_1=\nabla f /|\nabla f |$ is an eigenvector of
$Rc$;

(3) the second fundamental form $h_{ab}$ of $\Sigma_c$ is of the
form $h_{ab}=\frac{H}{n-1} g_{ab}$;

(4) the mean curvature $H=\frac{(n-1)R_\rho}{|\nabla f|}$ is
constant on $\Sigma_c$;

(5) on $\Sigma_c$, the Ricci tensor of $(M^n, g_{ij}, f)$ either has
a unique eigenvalue $\lambda$, or has two distinct eigenvalues
$\lambda$ and $\mu$ of multiplicity $1$ and $n-1$ respectively. In
either case, $e_1=\nabla f /|\nabla f |$ is an eigenvector of
$\lambda$. Moreover, both $\lambda$ and $\mu$ are constant on
$\Sigma_c$.

}

\proof Under this chosen orthonormal frame, we have $f_1=|\nabla f|$
and $f_2=f_3=\cdots=f_n=0$.  When $D_{ijk}=0$, we have from
Proposition 2.3 that
\begin{equation}\label{Proof22}
R_{1a}=0
\end{equation}
and
\begin{equation}\label{Proof23}
R_{ab}=\frac{R-R_{11}}{n-1}g_{ab}.
\end{equation}
Therefore, we obtain (1) from applying \eqref{Proof22} to
\eqref{Proof4} and \eqref{Proof5}, respectively. In particular, (2)
can be obtained from \eqref{Proof22} directly.

By the definition of $h_{ab}$, we have
\begin{equation}\label{Proof24}
h_{ab}=\langle\nabla_{e_a}\Big(\frac{\nabla f}{|\nabla
f|}\Big),e_b\rangle =\frac{1}{|\nabla
f|}f_{ab}=\frac{R_\rho}{|\nabla f|}g_{ab},\end{equation} where the
last equality comes from \eqref{Proof2}. Hence,
$H= h_{ab}g^{ab}=\frac{(n-1)R_\rho}{|\nabla f|}$ is constant from both $R$
and $|\nabla f|$ constant on $\Sigma_c$. Thus, (3) and (4) are
proved.

Using \eqref{Proof5}, we have
\begin{equation}\label{Proof25}
R_{1}=\frac{1}{m}R_\rho f_1-\frac{1}{n-1}R_{11}f_1=|\nabla
f|\Big(\frac{1}{m}R_\rho -\frac{1}{n-1}R_{11}\Big).
\end{equation}
From the definition of covariant derivative, we have
\begin{equation}\label{Proof26}\aligned
R_{,1a}=&e_1e_a(R)-\nabla_{e_a}e_1(R)\\
=&e_1e_a(R) -\langle\nabla_{e_a}e_1,e_1\rangle
R_1-\langle\nabla_{e_a}e_1,e_b\rangle R_b\\
=&e_1e_a(R)\\
=&0\endaligned\end{equation} since $R$ is constant on $\Sigma_c$.
Hence, we obtain from \eqref{Proof25}
\begin{equation}\label{Proof27}
0=R_{,1a}=-\frac{|\nabla f|}{n-1}R_{11,a}.
\end{equation}
Applying
$$R_{11,a}=e_a(R_{11})-2R(\nabla_{e_a}e_1,e_1)
=e_a(R_{11})-2h_{ab}R_{1b}=e_a(R_{11})$$ to \eqref{Proof27} yields
\begin{equation}\label{Proof28}
e_a(R_{11})=0,
\end{equation} which shows that $\lambda=R_{11}$ is constant on
$\Sigma_c$. By means of \eqref{Proof23} we know that for distinct
$a$, the eigenvalues of $R_{aa}$ are the same. Hence, we have the
eigenvalue $\mu$ is also constant. This completes the proof of
Proposition 2.4.

\noindent{\bf Proposition 2.5.} {\it Let $(M^n,g,f)$ be a quasi
Yamabe gradient soliton satisfying \eqref{Int3}, and let
$\Sigma_c=\{x|f(x)=c\}$ be the level surface with respect to regular
value $c$ of $f$.

(1) If $n=3$, then the sectional curvature of $\Sigma_c$ with the
induced metric is constant;

(2) If $n=4$ and $D_{ijk}=0$, then $W_{ijkl}=0$ on $\Sigma_c$.
Moreover, the sectional curvature of $\Sigma_c$ with the induced
metric is constant.

(3) If $n\geq 5$ and $W_{ijkl}=0$, then the sectional curvature of
$\Sigma_c$ with the induced metric is constant.

}

\proof It is well known that, for $n=3$, the Weyl curvature tensor
$W_{ijkl}$ vanishes identically. Hence, we obtain $D_{ijk}=0$ from
Proposition 2.2. Under the chosen local orthonormal frame as in
Proposition 2.4, we have from the Gauss equation and \eqref{Proof7}
and \eqref{Proof24}, for $a\neq b$:
\begin{equation}\label{Proof29}\aligned
R^{\Sigma}_{abab}=&R_{abab}+h_{aa}h_{bb}-h_{ab}^2\\
=&R_{abab}+\frac{(R_\rho)^2}{|\nabla f|^2}\\
=&R_{aa}+R_{bb}-\frac{R}{2}+\frac{(R_\rho)^2}{|\nabla f|^2}\\
=&\frac{R}{2}-R_{11}+\frac{(R_\rho)^2}{|\nabla f|^2}
\endaligned\end{equation} is constant. It completes the proof of (1).

Next, we prove (2). Since $D_{ijk}=0$, we have from
Proposition 2.2 $$W_{ijkl}f_l=0.$$ Hence, on the level surface
$\Sigma_c$, we have
\begin{equation}\label{Proof30} W_{ijk1}=0,
\ \ \  {\rm for}\ 1\leq i,j,k\leq 4.\end{equation} It remains to
show that
\begin{equation}\label{Proof31} W_{abcd}=0,
\ \ \  {\rm for}\ 2\leq a,b,c,d\leq 4.\end{equation} This
essentially reduces to showing the Weyl curvature tensor is equal to
zero in 3 dimensions (see \cite{Hamilton82}, p.276-277 or
\cite{CaoChen1}, p.13). Therefore, we have $W_{ijkl}=0$. When $n\geq
4$, the proof for constant sectional curvature is similar. We omit
it here. It concludes the proof of Proposition 2.5.

\noindent{\bf Proof of Theorem 1.1.} Following the proof of
Daskalopoulos and Sesum in \cite{Sesum11}, on the level surface
$\Sigma_c=\{x|f(x)=c\}$, we can express the metric $ds^2$ as
$$ds^2=\frac{1}{|\nabla f|^2}(f,\theta)df^2+g_{ab}(f,\theta)d\theta^ad\theta^b,$$
where $\theta=(\theta^2,\cdots,\theta^n)$ denotes intrinsic
coordinates for $\Sigma_c$. From \eqref{Proof29}, we know that
$\Sigma_c$ also has positive sectional curvature with respect to the
induced metric. Moreover, Proposition 2.4 and Proposition 2.5 show
that $\frac{1}{|\nabla f|^2}(f,\theta)=\frac{1}{|\nabla f|^2}(f)$
and $g_{ab}(f,\theta)=g_{ab}(f)$. Let $S$ be the set of critical
points of $f$. Then the measure of $S$ is zero. Hence, on
$M^n\setminus S$, we have $$ds^2=\frac{1}{|\nabla
f|^2}(f)df^2+g_{ab}(f)d\theta^ad\theta^b,$$ which shows that
$(M^n,g,f)$ is rotationally symmetric by using the arguments as in
\cite{Sesum11}. We complete the proof of Theorem 1.1.

\section{Proof of Theorem 1.2}

For any smooth function $u$ on $(M^n, g)$, we introduce the
following linear differential operator
\begin{equation}\label{oper1}
L(u)={\rm div}_{m,f}(\nabla u):=e^{\frac{f}{m}}{\rm
div}(e^{-\frac{f}{m}}\nabla u).
\end{equation} Then we have the following:

\noindent{\bf Lemma 3.1.} {\it  Let $(M^n, g)$ be a compact
Riemannian manifold. Then we have
\begin{equation}\label{oper2}
\int\limits_{M^n}v L(u)\,d\mu=\int\limits_{M^n}u L(v)\,d\mu,\ \ \ \
\forall\, u,v\in \,C^\infty(M^n)\end{equation} where $d\mu$ denotes
the measure $e^{-\frac{f}{m}}\,dV_g$. That is, the linear
differential operator $L$ is self-adjoint with respect to $L^2$
inner product under the measure $d\mu$. In particular, for any
smooth function $u$, we have
\begin{equation}\label{oper3}
\int\limits_{M^n}L(u)\,d\mu=0.\end{equation} }

\proof A direct calculation gives \begin{equation}\label{addaddadd}\aligned \int\limits_{M^n}v
L(u)\,d\mu=-\int\limits_{M^n}v_iu_i\,d\mu=\int\limits_{M^n}u
L(v)\,d\mu.
\endaligned\end{equation} This shows
that $L$ is self-adjoint with respect to $L^2$ inner product under
the measure $d\mu$. \eqref{oper3} is a special case when $v=1$ in
\eqref{oper2}. We complete the proof of Lemma 3.1.

Now we come back to prove Theorem 1.2. From \eqref{oper1} and  \eqref{Proof3},  $L(f)=\Delta
f-\frac{1}{m}|\nabla f|^2=nR_\rho$, thus we have
$$n\int\limits_{M^n}R_\rho\,d\mu=\int\limits_{M^n}L(f)\,d\mu=0,$$
which shows that
\begin{equation}\label{oper4}
\int\limits_{M^n}R_\rho\,d\mu=0.\end{equation} By Lemma 2.1, we have by using the second Bianchi identity and by choosing the local orthogonal frame $\{e_1,\cdots,e_n\}$,
$$\aligned \Delta(R_\rho)=&\frac{1}{m}(R_\rho)_i
f_i+\frac{1}{m}R_\rho
f_{ii}-\frac{1}{n-1}R_{ij,i}f_j-\frac{1}{n-1}R_{ij}f_{ij}\\
=&\frac{1}{m}(R_\rho)_i f_i+\frac{1}{m}R_\rho(nR_\rho
+\frac{1}{m}|\nabla
f|^2)-\frac{1}{2(n-1)}R_{j}f_j-\frac{1}{n-1}R_{ij}f_{ij},
\endaligned$$ which gives
\begin{equation}\label{2Proof3}
\aligned L(R_\rho)=&\Delta(R_\rho)-\frac{1}{m}(R_\rho)_i
f_i\\
=&\frac{n}{m}(R_\rho)^2+\frac{1}{m^2}|\nabla
f|^2R_\rho-\frac{1}{2(n-1)}(R_\rho)_i f_i-\frac{1}{n-1}R_{ij}f_{ij}.
\endaligned\end{equation}
Applying Lemma 2.1 again, we have
\begin{equation}\label{2Proof4}\aligned
-\frac{1}{n-1}R_{ij}f_{ij}=&-\frac{1}{n-1}R_{ij}(R_\rho\,g_{ij}
+\frac{1}{m}f_if_j)\\
=&-\frac{1}{n-1}R_\rho\,R+\frac{1}{m}f_i[(R_\rho)_i-\frac{1}{m}R_\rho\,f_i]\\
=&-\frac{1}{n-1}R_\rho\,R+\frac{1}{m}(R_\rho)_if_i-\frac{1}{m^2}|\nabla
f|^2R_\rho.
\endaligned\end{equation}
Inserting \eqref{2Proof4} into \eqref{2Proof3}, we obtain
\begin{equation}\label{2Proof5}\aligned
L(R_\rho)=&[\frac{1}{m}-\frac{1}{2(n-1)}](R_\rho)_i
f_i+\frac{n}{m}(R_\rho)^2-\frac{1}{n-1}R_\rho\,R.
\endaligned\end{equation}
Integrating \eqref{2Proof5} on $M^n$, we obtain by use of \eqref{Proof1}, \eqref{addaddadd} and \eqref{oper4}
\begin{equation}\label{2Proof6}\aligned
0=&\int\limits_{M^n}L(R_\rho)\,d\mu\\
=&\int\limits_{M^n} \Big\{[\frac{1}{m}-\frac{1}{2(n-1)}](R_\rho)_i
f_i+\frac{n}{m}(R_\rho)^2-\frac{1}{n-1}R_\rho\,R\Big\}\,d\mu\\
=&[\frac{1}{m}-\frac{1}{2(n-1)}]\int\limits_{M^n}(R_\rho)_i
f_i\,d\mu+[\frac{n}{m}-\frac{1}{n-1}]\int\limits_{M^n}(R_\rho)^2 \,d\mu\\
=&-[\frac{1}{m}-\frac{1}{2(n-1)}]\int\limits_{M^n}R_\rho L(
f)\,d\mu+[\frac{n}{m}-\frac{1}{n-1}]\int\limits_{M^n}(R_\rho)^2
\,d\mu\\
=&\frac{n-2}{2(n-1)}\int\limits_{M^n}(R_\rho)^2 \,d\mu,
\endaligned\end{equation} where we used $L(f)=nR_\rho$ in the last equality.
Clearly, \eqref{2Proof6} shows that $R_\rho=0$. Hence, we obtain
from \eqref{Proof3}
\begin{equation}\label{2Proof7}
\Delta f-\frac{1}{m}|\nabla f|^2=0.
\end{equation} Thus we get that $f$ must be constant since $M^n$ is compact.
We complete the proof of Theorem 1.2.

\noindent {Guangyue Huang} \\
\noindent Department of Mathematical Sciences, Tsinghua University,
Beijing 100084, P.R. China. \ \  E-mail
address:~\textsf{gyhuang$@$math.tsinghua.edu.cn }

\vspace*{1mm}
\noindent {Haizhong Li} \\
\noindent Department of Mathematical Sciences, Tsinghua University,
Beijing 100084, P.R. China. \ \  E-mail
address:~\textsf{hli$@$math.tsinghua.edu.cn }

\end{document}